\documentclass[11pt,titlepage]{amsart}
\usepackage{fancybox, fancyvrb}
\usepackage[dvipsnames]{xcolor}
\usepackage{  relsize}
\usepackage{graphicx}
\usepackage{amssymb, amsmath}
\usepackage{euler}
\usepackage{fancybox, fancyvrb}

\renewcommand{\theenumi}{\Alph{enumi}}

\newcommand{\thalf}{\tfrac12}
\font\small=cmr8 scaled \magstep0
\font\smallish=cmr8 scaled \magstep1

\outer\def\beginsection#1\par{\medbreak\bigskip
      \message{#1}\leftline{\bf#1}\nobreak\medskip
\vskip-\parskip
      \noindent}
\newcommand{\divisors}{\rm{div}}
   
\def\laq{\raise 0.4ex\hbox{$<$}\kern -0.8em\lower 0.62
ex\hbox{$\sim$}}
\def\gaq{\raise 0.4ex\hbox{$>$}\kern -0.7em\lower 0.62
ex\hbox{$\sim$}}
\def\beq{\begin{equation}}
\def\eeq{\end{equation}}
\def\bea{\begin{eqnarray}}
\def\eea{\end{eqnarray}}
\def\bean{\begin{eqnarray*}}
\def\eean{\end{eqnarray*}}

\DefineVerbatimEnvironment{code}{Verbatim}{fontsize=\relsize{-2}, frame=lines,%
fontfamily=courier,%
framesep=2\fboxsep,%
framerule=.5mm,%
rulecolor=\color{NavyBlue},%
fontseries=b,%
formatcom=\color{BlueViolet}}

\makeatletter
 \def\alphenumi{%
  \def\theenumi{\alph{enumi}}%
  \def\p@enumi{\theenumi}%
  \def\labelenumi{(\@alph\c@enumi)}}
\makeatother
\renewcommand{\labelenumi}{(\roman{enumi})}

\def \b {\beta}
\def \a {\alpha}

\def \b {\beta}
\def \a {\alpha}

\newcommand{\bs}{\boldsymbol}

\newcommand{\diff}{{\mathrm d}}

\usepackage {enumitem}

\DefineVerbatimEnvironment{code}{Verbatim}{fontsize=\relsize{-2}, frame=lines,%
  fontfamily=courier,%
  framesep=2\fboxsep,%
  framerule=.5mm,%
  rulecolor=\color{NavyBlue},%
  fontseries=b,%
  formatcom=\color{BlueViolet}}

\begin{document}
\author[  Enrico Onofri, 
  I.N.F.N., Gruppo Collegato di Parma]{  Enrico Onofri\\
  I.N.F.N., Gruppo Collegato di Parma, Parma, I-43124\\
\small   Department of Mathematics, Physics and Computer Science,
 Universit\`a di Parma, Italy
}

\title[Efficient Legendre transforms]{ Efficient Legendre polynomials transforms:\\
  from recurrence relations to Schoenberg's theorem}

\address{Dipartimento di Scienze Matematiche, Fisiche e Informatiche,
 Universit\`a di Parma, Parco Area delle Scienze 7,  43124 Parma, Italy}
\urladdr{www.pr.infn.it/~enrico.onofri}
\email{enrico.onofri@pr.infn.it}

\begin{abstract}
We report results on various techniques which allow to
  compute the expansion into Legendre (or in general Gegenbauer)
  polynomials in an efficient way. We describe in some detail the
  algebraic/symbolic approach already presented in Ref.\cite{VYO2017}
  and expand on an alternative approach based on a theorem of
  Schoenberg \cite{Sch1942}.
\end{abstract}

\maketitle

\section{Introduction}
This investigation stems from a problem regarding unitarity of the
pion-pion scattering amplitude in a model proposed by Veneziano and
Yankielowicz\footnote{See \cite{VYO2017} for motivation and a thorough
description of the model}: one is given the function
\begin{equation}
\label{eq:def}
A(x,M,\a) \equiv \sum_{k\in\divisors{(M\!+\!1)}}\,c_k(\beta,\gamma)\,
\frac{1-\a\!+\!M+t}{q_k!}\;
\frac{\Gamma(q_k+(\a+t)/k)}{\Gamma((\a+t)/k)}
\end{equation}
where the sum runs over the divisors of $M\!\!+\!\!1$,
$t=\thalf(x-1)(1-\a+M)$, $x=\cos\vartheta$ ($\vartheta$ is the
scattering angle), $\a\in(0,1)$ is the Regge slope, $q_k=(M+1)/k-1$
and $M=s+\a-1$ is a non--negative integer\footnote{$t$ and $s$ are
  Mandelstam's scattering parameters, but here $x$ and $M$ are used
  instead.}. The coefficients $c_k$ are assumed to be given by
$c_k\propto k^{-(1+\beta)}\,(1+(\log k)^{-\gamma})$, $\beta$ is a real
positive parameter ($\beta =2$ is the most popular value) and usually
one assumes $\gamma=\beta\!+\!1$.  The problem is to determine the
range of the real parameters $(\a,\beta,\gamma)$ for which the 
coefficients of the expansion
$${\mathcal A}(s,t)\equiv A(x,M,\a) = \sum_{n\ge 0} a_n\,P_n(x)$$ into 
Legendre polynomials are all non--negative.

At first sight the problem looks  rather trivial; one has only to
compute the coefficients by Legendre--Fourier inversion
\begin{equation}
  \label{eq:1}
  a_n \propto \int_{-1}^1\,\diff{ x} \; A(x,M,\a) \,P_n(x)\;.
\end{equation}
and check their sign.  However, in practice, one has to face a variety
of problems. The number of integrals to be computed is very high if
one wants to explore the ``critical landscape'' parametrized by the
pairs $(\beta, M)$ and identify the critical value
$\a_{\rm crit}=\sup_M\a(\beta,M)$ such that positivity occurs for
$\a > \a_{\rm crit}$ and for all $M$. Moreover the function to be
analyzed is a high--degree polynomial which makes the integrals rather
sensitive to accuracy. Hence it is advisable to employ exact
arithmetic or at least extended precision. This is easily achieved
using a symbolic language like Mathematica or adopting an extended
precision package in Matlab.  We found that an alternative technique,
not requiring the calculation of even a single integral, is however
preferable if one wants to reach the result in a reasonable
time. Another strategy, which also does not require to compute Fourier
coefficients explicitly, has been developed based on a theorem of Schoenberg
which relates the non--negativity of the expansion coefficients to a
positivity property of some matrices built out of the function
$A(s,t)$. In the following sections we present the ``algebraic
approach'' and the ``Schoenberg's theorem'' approach. 
\section{The algebraic approach}
We note that $A(x, M, \a)$ is a polynomial of degree $M\!+\!1$ which is defined
as a sum of terms which are easily factorized in simple monomials,
being a quotient of gamma-functions (Pochammer's symbols). It follows
that an easy technique to derive its expansion into Legendre
polynomials is provided by the recurrence relation
\begin{equation}
  x\,P_n(x) = \frac{n\!+\!1}{2n\!+\!1}\,P_{n+1}(x) + \frac{n}{2n\!+\!1}\,P_{n-1}(x)\;.
\end{equation}
On the $d-$dimensional sphere $S^d$ one has to consider the recurrence
for Gegenbauer's polynomials, namely
\begin{equation}
  \label{eq:2}
   x\,C^\lambda_n(x) = \frac{n\!+\!1}{2(n\!+\lambda)}\,C^\lambda_{n\!+\!1}(x) \!+\!
\frac{n\!+2\lambda-1}{2(n\!+\lambda)}\,C^\lambda_{n-1}(x),\;\left( \lambda=\thalf(d-1)\right)\;.
\end{equation}
Let $A_k=\prod_j(x-x_j)$ be the
factorization of the $k-th$ term in Eq.\eqref{eq:def}, suitably
normalized; starting from $A_{k,0}(x)\equiv 1$ and applying the first
factor we get
\begin{equation*}
  A_{k,1} = (x-x_1) = P_1 - x_1\,P_0
\end{equation*}
and, subsequently
\begin{eqnarray*}
  A_{k,2} &= (x-x_2)\,A_{k,1} = x(P_1 - x_1\,P_0)-x_2 (P_1 -
  x_1\,P_0)\\
&= \frac23 P_2  -(x_1+x_2) P_1 + (x_1 x_2+\frac13) P_0
\end{eqnarray*}
and so forth,  so that after applying this formula a number of times
to reach the degree of the polynomial, we get by inspection the required
expansion.  The process is easily realized by using a symbolic
language like Mathematica or {\tt form}; alternatively we can
represent the multiplication by a factor $x$ as an infinite matrix
\begin{equation*}\label{eq:X}
  {\mathcal X} =
  \begin{pmatrix}
    0&\frac13&0&0&0&0&\ldots\\
    1&0&\frac25&0&0&0&\ldots\\
    0&\frac23&0&\frac37&0&0&\ldots\\
    0&0&\frac35&0&\frac49&0&\ldots\\
    \vdots&&&\ddots&&\ddots&
  \end{pmatrix}
\end{equation*}
and in general for any dimension, using Eq.~\eqref{eq:2}.
The expansion is then computed by applying the product
$\prod_j({\mathcal X} - x_j)$ to the initial vector
$A_{k,0}=\{1,0,0,\ldots\}'$. Let us note that using a symbolic
approach makes it possible to derive the result by working with
symbolic parameters $(a,\beta)$ ; using a floating point
representation for the coefficients, instead, one loses this
possibility and moreover one has to pay attention to arithmetic
accuracy, since floating point errors can be amplified in the process.
Hence in a floating point approach one has to use multiple precision
arithmetic\footnote{We successfully used the package from  {\tt
advanpix.com}}.  Having developed an efficient algorithm to
compute the expansion coefficients, the range of values in the
$(\alpha,\beta)$ plane where positivity holds {\sl for any} $M$ can be
determined by the classical bisection method in the following way:

{\sl 
\begin{enumerate}[label=\alph*)]
\item fix a value for $M$ and $\beta$;
\item fix a tolerance $\varepsilon$, e.g. $\varepsilon=10^{-6}$;
\item start with an interval for $\alpha$, say $\alpha_{\rm min}=0$, $\alpha_{\rm max}=1$ such that
  positivity is known to hold for $\alpha=1$ and it does not hold for $\alpha=0$;
\item compute the coefficients at the midpoint
  $\a =\thalf(\a_{\rm min}+\a_{\rm max})$;
\item if all coefficients are positive, then set $\a_{\rm max}=\a$,
  otherwise set $\a_{\rm min}=\a$;
\item repeat until $|\a_{\rm max}-\a_{\rm min}|\le\varepsilon$ and keep
  the mean value as $\a_{\rm crit}(M,\beta)$ up to $O(\varepsilon)$;
\item loop on $M$ and $\beta$ .
\end{enumerate}
}

One is interested to identify the points $(\a,\beta)$ where positivity
holds for \emph{every} $M$, i.e. we have to compute
$\a_{\rm crit}(\b)=\max_M \a_{\rm crit}(M,\b)$.  This requires, of
course, a compromise, since we have to fix a maximum value of $M$ to
make the calculation in a finite time, but \emph{a posteriori} we see
that after some large $M\approx 100$ the profile stabilizes. Here we
borrow some pictures from Ref\cite{VYO2017}
(Fig.1 and 2).
\begin{figure}[b] 
\centering
\includegraphics[width=0.65\linewidth]{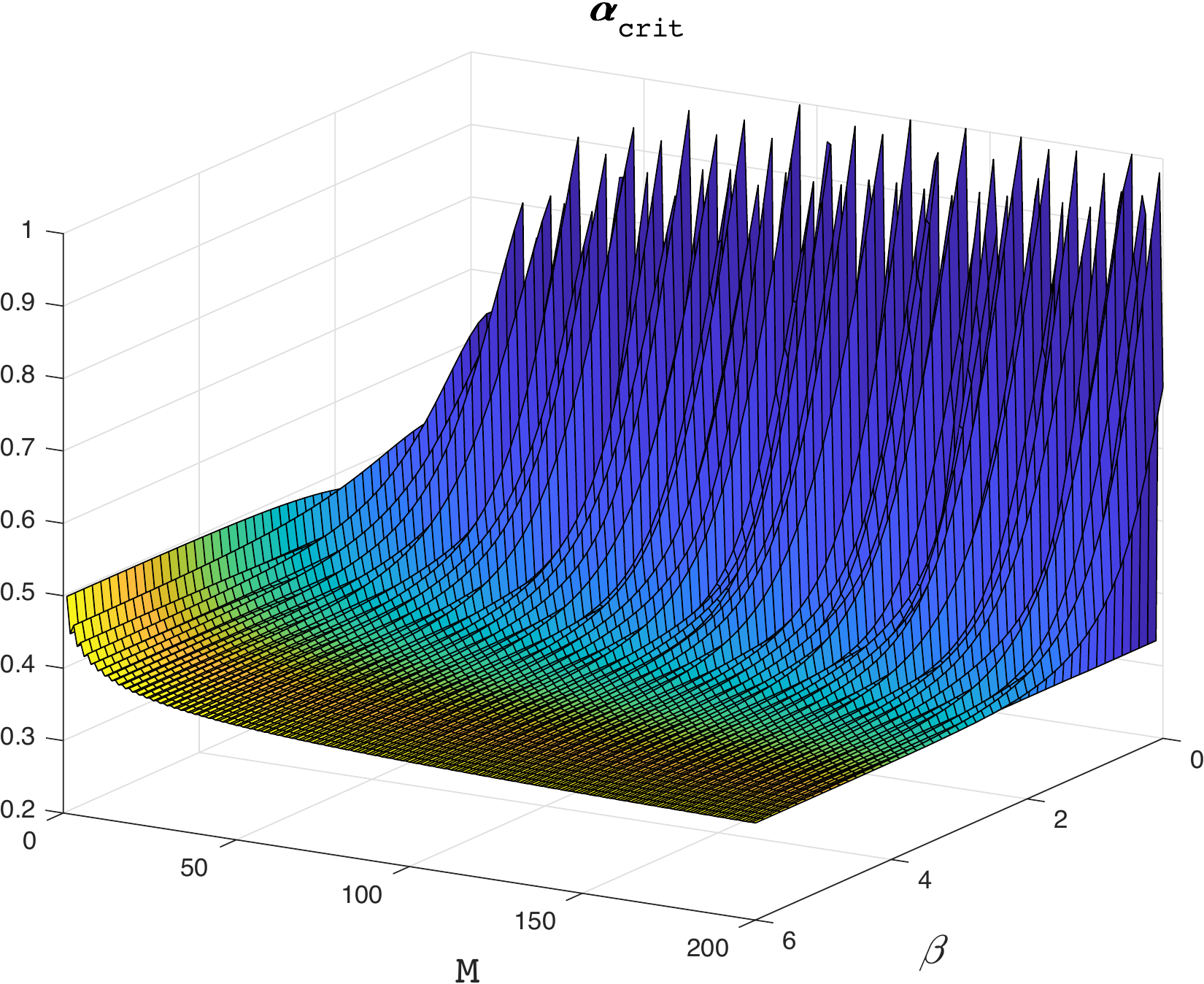}
\caption{{\smallish The ``critical landscape'': values of $\alpha$
    above the surface give positive residues (here and in the
    following {\footnotesize$\gamma=\beta+1$})}}
\label{fig:AcritAss}
\end{figure} 
\begin{figure}[b] 
\centering
  \includegraphics[width=.65\linewidth]{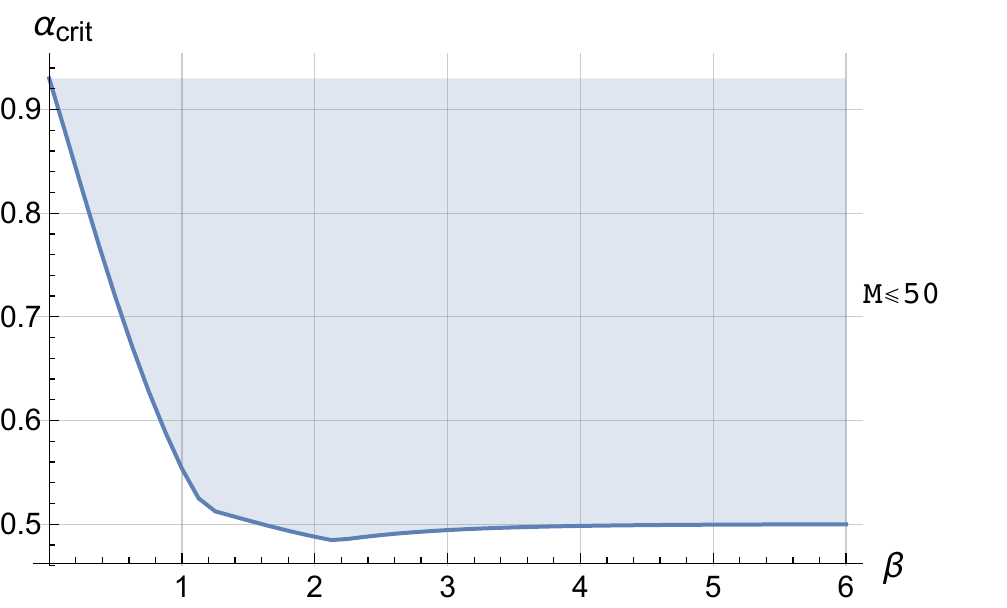}
  \caption{{\smallish Same as the plot of Fig.\ref{fig:AcritAss} projected on the
      ${\footnotesize (\beta,\alpha)}$ plane. Shaded region corresponds to positive
      residues.
}}
\label{fig:shaded}
\end{figure}

\subsection{Chebfun}
A very simple approach, which essentially makes it trivial to write
the floating point code, adopts the {\sl Chebfun\/} environment in
Matlab. Here the expansion coefficients are built in and we get the
result with essentially a single instruction. The bisection strategy
then is again used to get the critical values of $\a$. This option,
however, would require a multiprecision extension of Chebfun, which is
not available at present, and it does not appear to be around the
corner. Hence this method has been mainly used  for debugging purposes.

\subsection{Spectral method}
The matrix representing multiplication by $x$, namely $\mathcal X$
defined in Eq.\eqref{eq:X}, is known to have the zeros of the Legendre
polynomial $P_{N+1}(x)$ as eigenvalues and its eigenvectors are
exactly known in terms of Legendre polynomials of lower indices
evaluated at the eigenvalues\footnote{The result is due to Golub and
  Welsch \cite{GW1969}}. This means that we exactly know the spectral
representation of ${\mathcal X}$ and any function of it can be
computed in algebraic form. This is particularly useful if the
function to be expanded is not a polynomial. However in this case the
maximum order $N$ is not defined, hence it must be introduced as a
cutoff, which involves an error in the calculation of the coefficients
which must be estimated by extrapolating to large $N$. It turns out
that this spectral strategy coincides step by step with the idea of
computing the integrals in Eq.\eqref{eq:1} by Gaussian quadrature:
since this is builtin in Mathematica it does not cost much to try this
alternative strategy. Our experience is that we are not going to earn
much in terms of computing efficiency. See Ref.\cite{VYO2017} for
details.

\section{The alternative route: Schoenberg's theorem}
The main new contribution of this report regards a further strategy we may
adopt to solve the original problem. This is based on a theorem going
back to 1942 due to I.~J.~Schoenberg
\cite{Sch1942}\footnote{for a more recent review, see
  \cite{JStewart76}}: given a function $f(\cos\vartheta)$, the
positivity of the coefficients in its expansion into Legendre
functions is directly linked to a property of the function which is
know as ``positive definiteness'' (or ``semi-positive definitness'').
\begin{flushleft}
  {\bf Definition: }{\sl a function $f(\cos\vartheta)$ is said to be
    positive semi-definite on $S^d$ if for any integer $n>0$ and any
    $n$-tuple of unit vectors $\{\bs{v}_1,\ldots,\bs{v}_n\}$ on $S^d$
    the matrix $F_{ij}\equiv f(\bs{v_i}\cdot \bs{v}_j)$
    is positive semi-definite, i.e. $\sum_{ij} F_{ij}\,\a_i\,\a_j \ge  0$.\\[.5em]
    {\bf Theorem } (Schoenberg 1942): A function $f(\cos\vartheta)$ is
    positive semi-definite if and only if its expansion into Legendre
    functions has non-negative coefficients. }
\end{flushleft}
By ``{\sl Legendre functions\/}'' it is understood Legendre polynomials
$P_\ell(x)$ in dimension $d=2$, and Gegenbauer's polynomials
$C_n^\lambda$ in any dimension\footnote{notice that by ``dimension'' we always mean
  the dimension of the manifold $S^d$ embedded in $R^{d+1}$.}
$d$, with $\lambda=\thalf(d-1)$.  

The question then arises whether this equivalence can be exploited in
a practical way to prove positivity of the expansion
coefficients. There is an obvious difficulty: the theorem requires in
principle to test the positivity of the matrix $F_{ij}$ for
\emph{every} integer $n$ and this is clearly impossible in
practice. One is then led to ask whether it would be sufficient to
perform the test in a suitable range of values for $n$, say
$n<N_{\rm max}$.  Assuming that a negative coefficient is present at
some level $\ell$, is there a connection between $\ell$ and the range
of $n$ which will allow to detect the presence of this negative
coefficient?  We need  some condition similar to what is known
to hold in the theory of signals, where getting information about the
spectrum of a signal around a frequency $\nu$ requires sampling the
signal at least with a frequency $2\nu$ (Nyquist-Shannon theorem). In
the case of functions defined on the $d$-dimensional sphere what
should be the minimal dimension of the $n$-tuples of vectors which can
detect a negative component in the Legendre expansion? To investigate
this problem we performed a series of tests which point to a
preliminary result: ``{\sl if the first negative component
  is at level $\ell$ then one should check positivity with $n$-tuples
  such that $n$ is at least twice as large as the number of spherical
  harmonics of index less or equal to $\ell$}''.  In dimension $2$
this requires testing positivity of the matrix $F$ for $n$-tuples with
$n\ge 2(\ell+1)^2$, in higher dimensions the correct formula can be
found in Ref.\cite{Miller66} or in the book by Hochstadt
Ref.\cite{Hochstadt71} (see next section).

In a recent paper by {\sl Beatson, W. zu Castell and Xu} \cite{BzCXu2014} it is
conjectured that a necessary and suffcient condition for positivity for every
dimension is that the power expansion $f(x)=\sum c_k\,(\cos\vartheta)^k$ has all
positive Taylor coefficients. Unfortunately again we cannot take profit of this
conjecture, since the positivity condition \emph{for every $S^d$ (any $d$!)}
which is too strong a requirement to be applied to the amplitude given in Eq.(1)
and one readily realizes that the condition is \emph{not} respected for any
value of $\alpha$ signaling the fact that the amplitude fails positivity for
sufficiently high dimension (perhaps 25?).

Let us remark that the analysis in terms of Schoenberg's theorem
allows one to infer that, since the set of $n$-tuples in $d$ dimensions contains
all $n$-tuples in less than $d$ dimensions, positivity at a given value of the
parameters for a given function $f(x)$ in dimension $d$ {\sl implies positivity
  for the same function for all dimensions lower than $d$ }. This property has
been checked by the algebraic method but is far from obvious looking at the
recursion relations. On the other hand this property may be difficult to check
numerically, the problem being that going up with $d$ the detection of a
negative component requires to go to higher dimensional $n$tuples, typically
$n\propto \ell_0^d$. This means that our estimate of a critical $\alpha_c$ in
$d-$dimensions can be affected by a systematic error much higher than in lower
dimensions. On the other hand, if we choose the same function $f(x)$ and apply
the test in various dimensions, we check that positivity in dimension $d_>$
implies positivity in dimension $d$ with $d<d_>$.  From this point of view,
let's say that from an {\sl algorithmic\/} point of view, the
statement/conjecture of Ref.~\cite{BzCXu2014} is rather trivial: in the limit of
large dimension it is well-known that
$ \label{eq:limit} \lim_{d\to\infty}C_n^{(d-1)/2}(z)/d^n = z^n/n!\;, $ hence the
series of hyperspherical polynomials converges to a simple Taylor series (see
\cite{AAR99}, pag.303)!

\subsection{Numerical results}
The numerical test we have devised runs as follows (initially for
dimension 2)

\begin{enumerate}[label=\alph*)]
  {\sl
  \item build the function $f(x)$ directly as a finite sum of Legendre
    functions with known coefficients $c_\ell$ depending on a
    parameter $\alpha$, explicitly $c_\ell\propto 1-\a/\a_0$
    starting at some index $\ell_{\rm min}$;
  \item  extract  an $n$-tuple of unit vectors,
    $\{\bs{v}_1,\ldots,\bs{v}_n\}$ uniformly at random on $S^d$;
  \item evaluate the spectrum of the matrix
    $f(\bs{v}_i\cdot\bs{v}_j)$;
\item determine the value of $\a_0$ by a bisection strategy as
  discussed previously depending on the presence of negative
  eigenvalues;
\item for each value of $\ell_{\rm min}$ perform the test on a range
  $n\in\{1,2\,(\ell_{\rm min}\!+\!1)^2\}$ and check the accuracy in determining
  $\a_0$;
\item repeat (b) a number of times  and gather the statistics;
\item  examine the dependence of the estimate for
  $\a_0$ depending on $n\;.$
\item Extend the experiment to higher dimensions.
}
\end{enumerate}

\subsubsection{Test runs}
Here we report on a number of test runs in dimension 2, 3, and higher. In the
first three runs we examine the effect of applying Schoenberg's test to
functions characterised by the same coefficients in the expansion of
hyperspherical polynomials irrespective of dimension. We stress that in this
case the coefficients are the same but the functions are different. The critical
value for $\a$ is the same.  The pattern which emerges is the following: for
small sets of unit vectors (small $n$) the negative coefficients are often
overlooked and as a result the critical $\a_0$ is overestimated; as $n$ grows
the estimate is better and better. If the occurrence of negative components
starts at $\ell$ the value of $n$ must be compared to $\ell^{d}$; a rationale
for this is probably that the number $H(\ell,d)$ of (hyper)spherical harmonics
of index up to $\ell$ on $S^d$ grows as $O(\ell^d)$, precisely
\begin{equation*}
\begin{aligned}
  H(\ell,2)&=(\ell+1)^2\\
  H(\ell,3)&=(\ell+1)(2+\ell)(2\ell+3)/6\\
  H(\ell,4)&=(\ell+1)(\ell+2)^2(\ell+3)/12\\
  H(\ell,5)&= (1 + \ell]) (2 + \ell) (3 +\ell) (4 + \ell) (5 + 2 \ell)
  /120\\
  ...&...
\end{aligned} 
  \end{equation*}
 etc. \footnote{The general formula for the number of linearly
  independent spherical harmonics on $S^d$ of degree $\ell$ is given by
  $N(d,\ell)=\frac{2\ell+d-1}{\ell+d-1} {{\ell+d-1}\choose{d-1}}$ (see
  \cite{Hochstadt71}).}

\begin{figure}
  \includegraphics[width=0.75\linewidth]{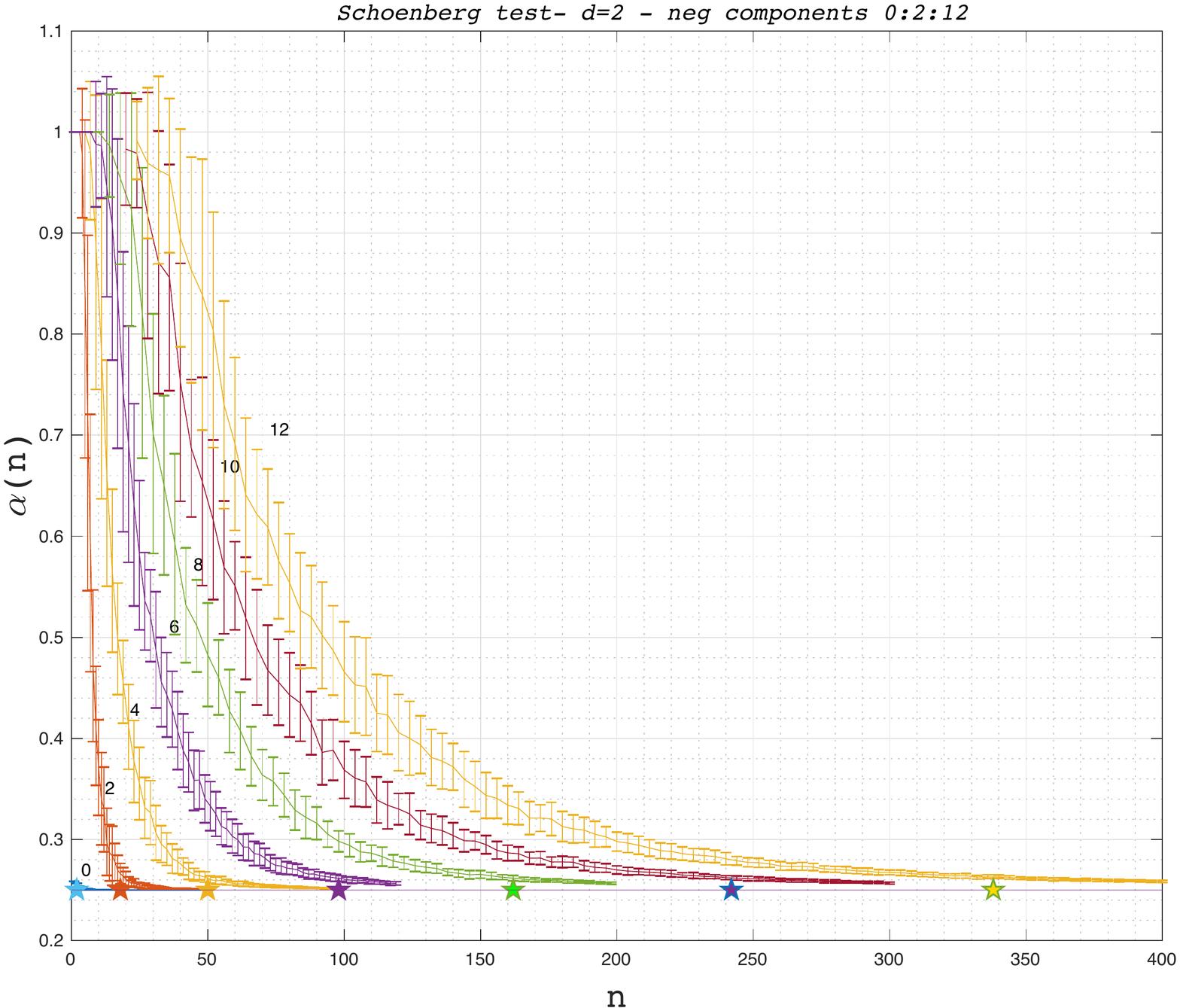}
  \caption{{\smallish Schoenberg's test at work: $\a_0=0.25$}}
\label{fig:SchoenbTest}
\end{figure} 

In this experiment we find that the critical value $\a_0$ is computed with an
error $\approx 5\%$ for $n\approx 2(\ell_{\rm min}+1)^2$; these values are shown
by the symbols $\star$ on the line at $\a=0.25$. It is clear that there is a
systematic error for small $n$ which dominates the statistical error represented
by the errorbars. In a real problem we do not know the value of
$\ell_{\rm min}$, hence one has to compute at the highest value reachable for
$n$ and one should try to extrapolate at infinite $n$.

In the following two plots we report the same test performed for dimensions
$d=3$ and $4$. The $x$-axis is scaled as $n^{1/3}$ and $n^{1/4}$, respectively,
to make the plot more compact; this is based on the fact that a good range for
$n$ starts at $2H(\ell_{\rm min},d) = O(\ell_{\rm min}^d)$.

In Fig.\ref{fig:dimen} we see many values of $d$ at the same time: also
here the function \emph{is not the same\/}, rather the coefficients of the
expansion in hyperspherical polynomials are the same, and we see how convergence
is slower in higher dimensions. To check Schoenberg test \emph{in all\/}
dimensions is obviously impossible, but here we see that the approximation
improves with the number of vectors in the $n$-tuple.

Finally in Fig.\ref{fig:dimen-incl} we report the result in the case
of the same function tested in different dimensions. As expected the
limiting values which approach the critical values of the parameter
$\a$ are ordered in such a way that positivity in $d$ implies
positivity in lower dimensions. This property applies to the limiting
values, while there is no reason to expect that it apply to any value
of the dimension of the $n$-tuples.  We applied the technique based on
Schoenberg's theorem to the original problem of finding the positivity
range for $\a_0$. A preliminary result is encouraging (see
Fig.\ref{fig:alpha}). However, if applicable, the algebraic method is
much more efficient. The alternative tecnique based on Schoenberg's
theorem may turn out to be the right choice in cases where an explicit
analytic expression is not available for the function $f(x)$.

\section*{Acknowledgments}

It's a pleasure to thank \emph{G. Veneziano} and \emph{S. Yankielowicz}
for initially suggesting the problem and for subsequent encouragment and advice. I thank
the Director of the SMFI Department, \emph{Prof. R. De Renzi},  for his generous hospitality.

\appendix
\section{The matlab code for  Schoenberg's test }
{\sl We report here relevant parts of the code. The complete working code can be
downloaded from the author's website\/}. 
\begin{code}[frame=lines, framesep=2\fboxsep, framesep=3mm,
 label='SchoenbTest.m'] 
  function [alpha, dalpha, setup, nout,Y] = SchoenbTest(setup)
  % Schoenberg test. A number of n-tuples are extracted for
  % each n in an interval fixed in setup. The mean and std over all runs 
  % is recorded. The output "alpha" contains the best
  % estimate of alpha as a function of the cardinality n of n-tuples and
  % dalpha gives the statistical error.  Y is a structure containing
  % the data needed by the various subroutines. The main loop on sweeps
  % is parallelized on the available cores (parfor)
  %
  % Usage:
  %       [alpha, dalpha, setup, nout , Y] = SchoenbTestS([setup])

if nargin<1,  setup = initialize();     % setup parameter
else
   setup = initialize(setup);     % setup parameter
end

Y.dim      = str2num(setup{1});
samples  = str2num(setup{2});  % putting in Y all parameters to be passed to subroutines
Y.tol        = str2num(setup{3});
NN           = str2num(setup{4}); % interval for n
Y.a0         = str2num(setup{5}); % critical value of alpha
Y.nmax   = str2num(setup{6});
hits          = str2num(setup{7});
cf             = str2num(setup{8}); % expansion coefficients

rng('shuffle','twister')        % initialize random numbers

Y.cf  = cf;
Y.cfP = zeros(size(cf));
Y.cfP(hits) = Y.cf(hits); % selects coefficients dependent on alpha

alpha  = zeros(length(NN),1);
dalpha = zeros(length(NN),1);
nout = zeros(length(NN),1);
index = 1;

for n = NN,
   acum = zeros(samples,1);
   
   parfor rep = 1:samples,                % parallelized loop
      dim = Y.dim(1);
      v = randn(dim,n);               % N vectors at random in R^{dim}
      s = ones(dim,1)*sum(v.^2);
      v = v./sqrt(s);
      z = v'*v;                         % matrix v_i . v_j
      a = bisec(z, Y, n);      
      acum(rep)=a;
   end

   nout(index) = n;    
   av = mean(acum);    dev = std(acum);
   alpha(index)   = av; dalpha(index)  = dev;
   index = index+1;
end

%------------------------------------------
function a = bisec(z, Y, n)

a1 = 0.0; a2 = 1.0;

while(a2-a1 > Y.tol)

   a = (a1+a2)/2;
   
   flag = negeig(a, z, Y, n);
   
   if flag,      a2 = a;
   else,         a1 = a;
   end

end


%%-------------------------------------------------
function flag = negeig(a, z, Y, N)
% sampling n-tuples on S^{dim-1} to study positivity of function F
% defined below. Fires a flag when a negative eigenvalue is detected:

flag = 0;
warning off

nu = Y.dim(2)/2-1;                       % parameter of Gegenbauer
                                          % polyn. sph.harm. on S^(dim-1)

C  = Y.cf;
P  = Y.cfP;
a0 = Y.a0;

C0 = 1;              % Gegenbauer computed by recurrence
C1 = 2*nu*z;
F  = C0*(C(1)-a/a0*P(1))*ones(size(z)); % P_0

for k = 1:Y.nmax,
   F = F + C1*(C(k)-a/a0*P(k));		% P_{k+1}(z) ....
   C2 = (2*(k+nu)*z.*C1 - (k-1+2*nu)*C0)/(k+1);
<   C0 = C1; 
   C1 = C2;
end

F = F.*(ones(size(F)));
F = 0.5*(F+F');

if(N<500),
   ev  = min(eig(F));
else
   opts.issym = true;
   opts.p = 40;
   ev = eigs(F, 1, 'sr', opts);
end

if(ev < 0),   
   flag = 1;  
end  % set flag and return

%%----- setting up parameters via "inputdlg"
function data = initialize(setup)   
% ......  instructions to set the parameters values interactively
data = inputdlg(entries, title, lines, default, AddOpts);
\end{code}
\bibliographystyle {amsalpha}

\begin{figure}[b] \label{Schoenb_II}
\centering
  \includegraphics[width=0.65\linewidth]{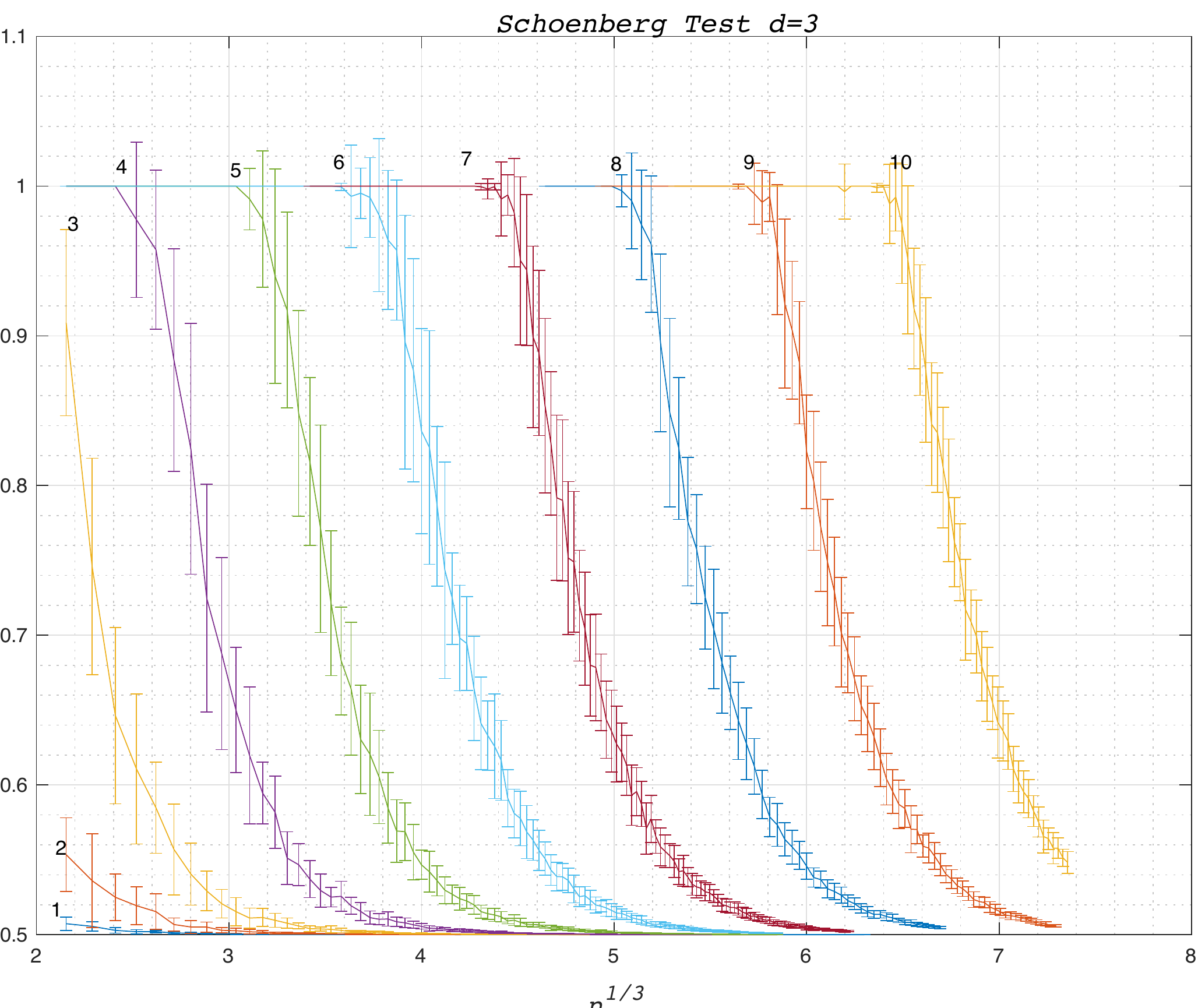}
  \caption{{\smallish Schoenberg's test at work: $\a_0=0.5$}}
\end{figure} 

\begin{figure}[b] \label{Schoenb_III} \centering
  \includegraphics[width=0.65\linewidth]{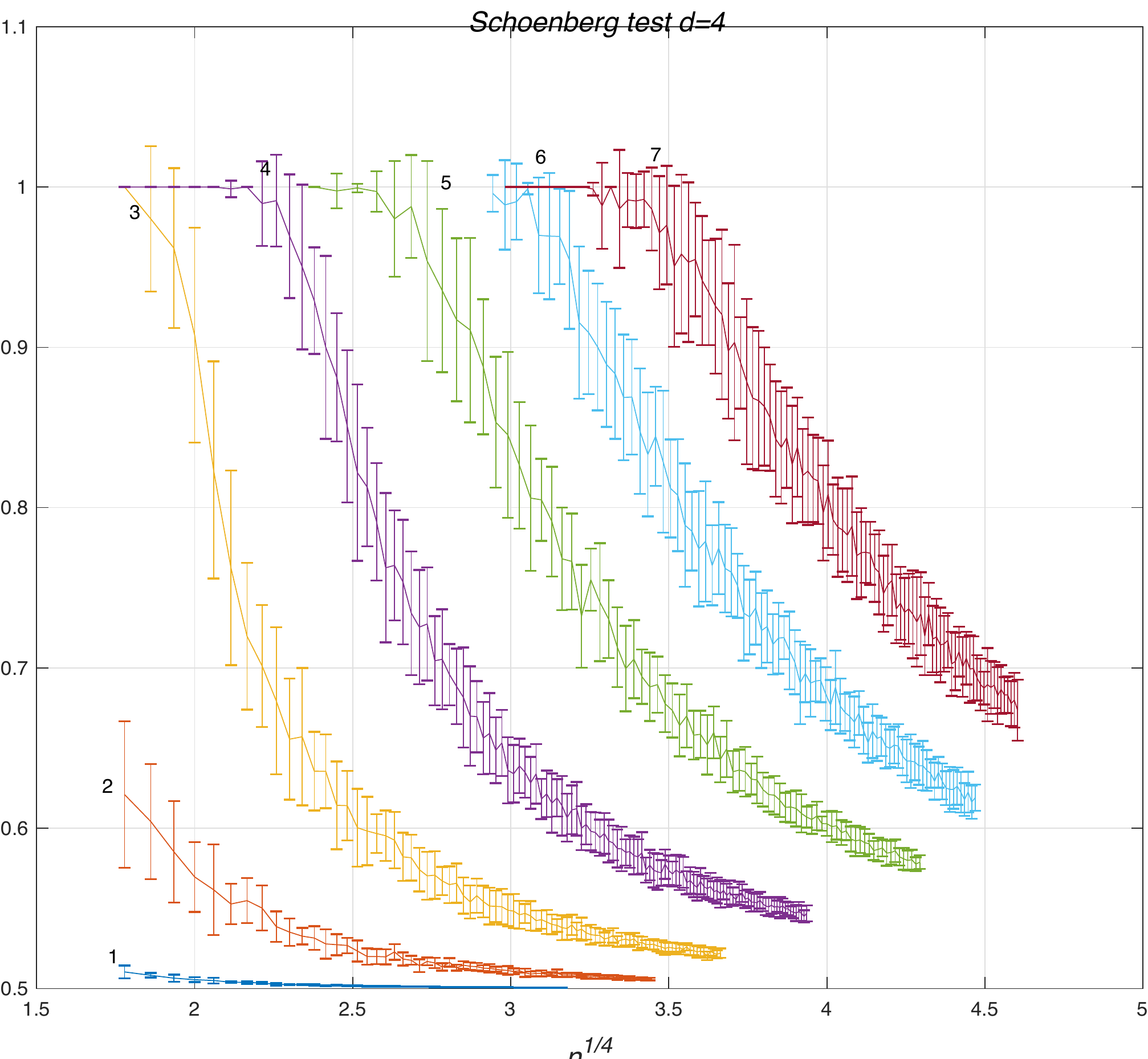}
  \caption{{\smallish Schoenberg's test at work: $\a_0=0.5$}}
\end{figure}

\begin{figure}[b]
\centering
  \includegraphics[width=0.65\linewidth]{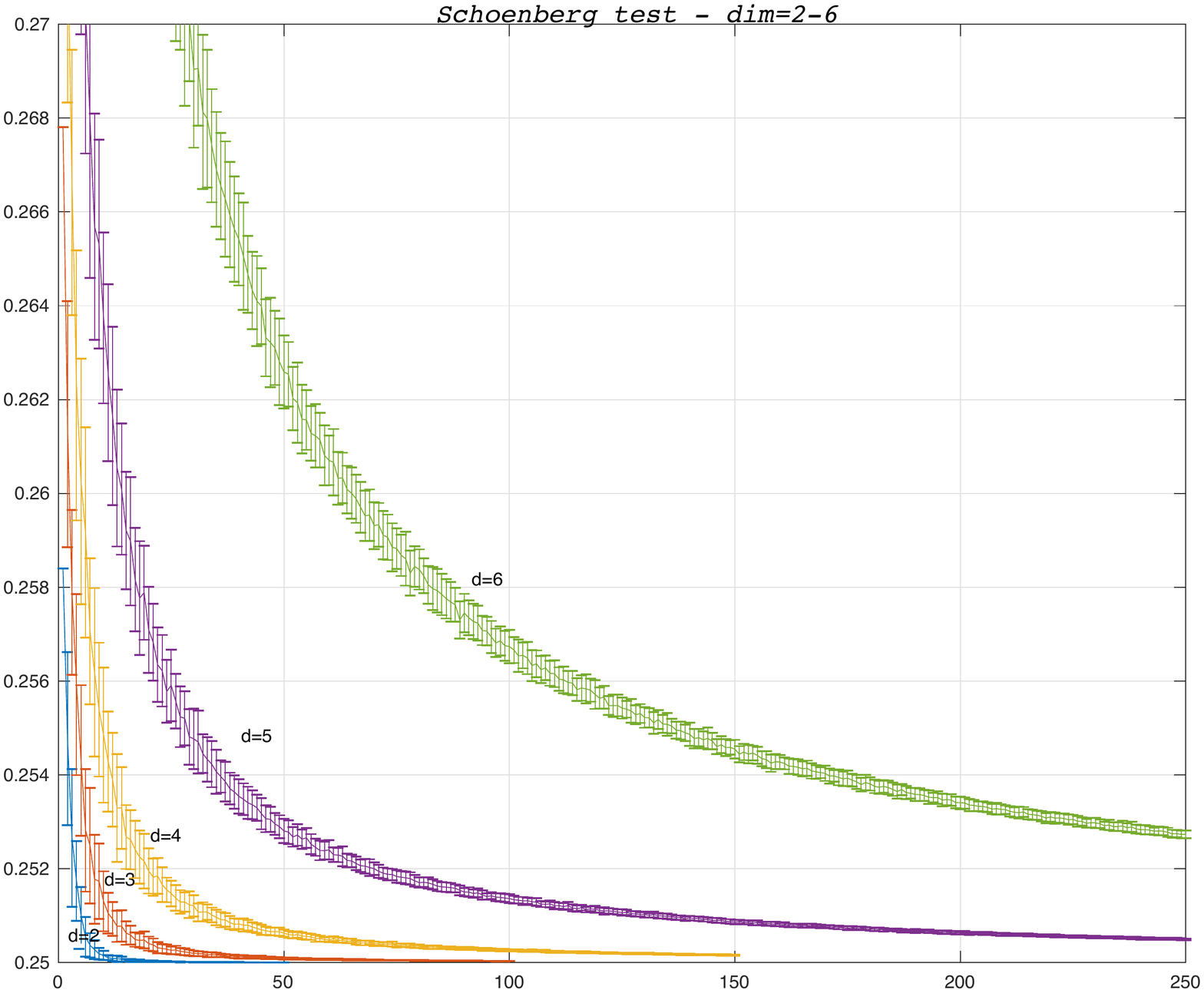}
\caption{{\smallish Exploring $n$-tuples in various dimensions at
    fixed parameters}. $\alpha_c=0.25$.}
\label{fig:dimen}
\end{figure} 
\begin{figure}[b]
\centering
  \includegraphics[width=0.8\linewidth]{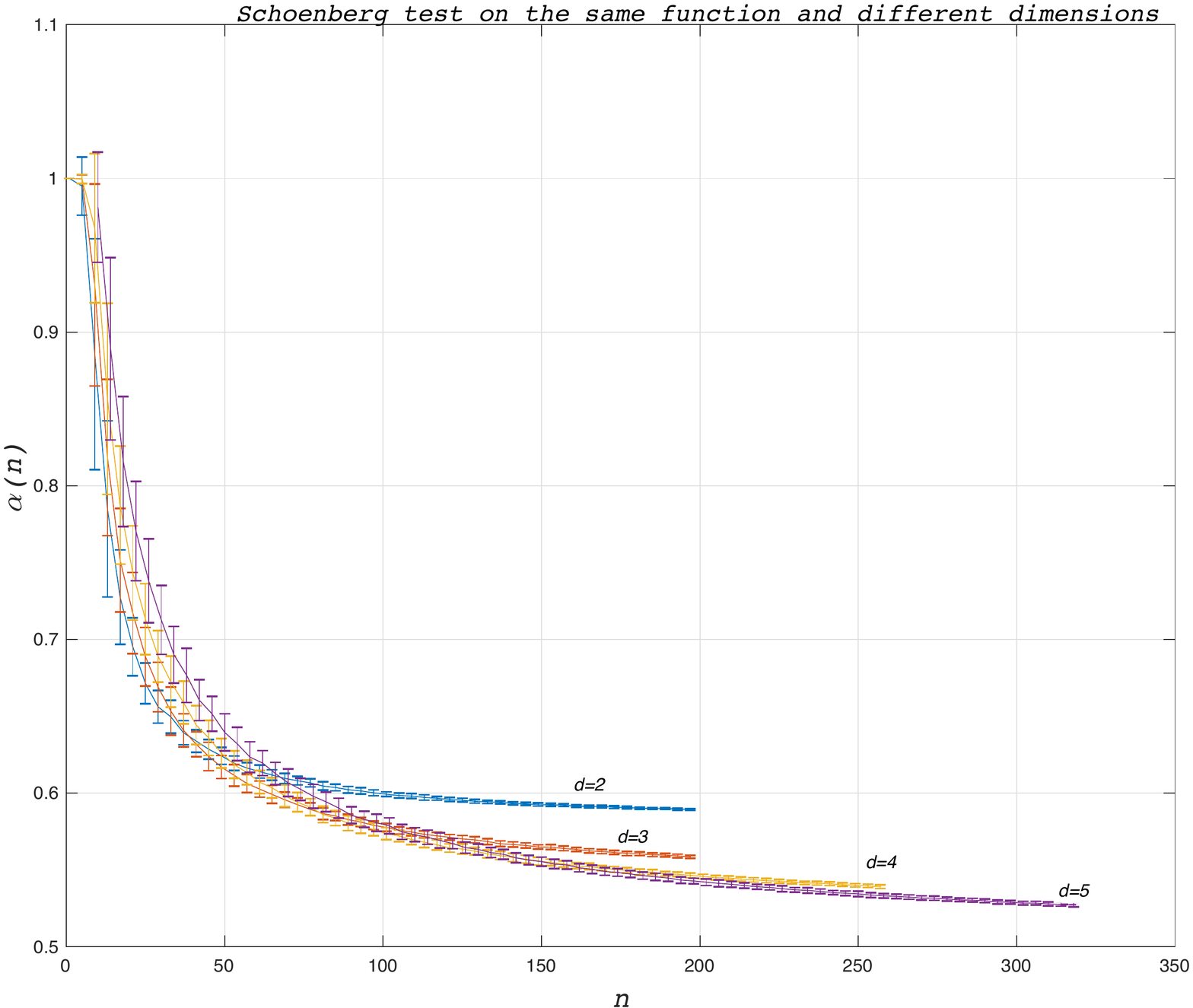}
  \caption{{\smallish Exploring $n$-tuples in various dimensions keeping the function $f(x)$
      fixed}. }
\label{fig:dimen-incl}
\end{figure} 
\begin{figure}[b]
\centering
\includegraphics[width=0.8\linewidth]{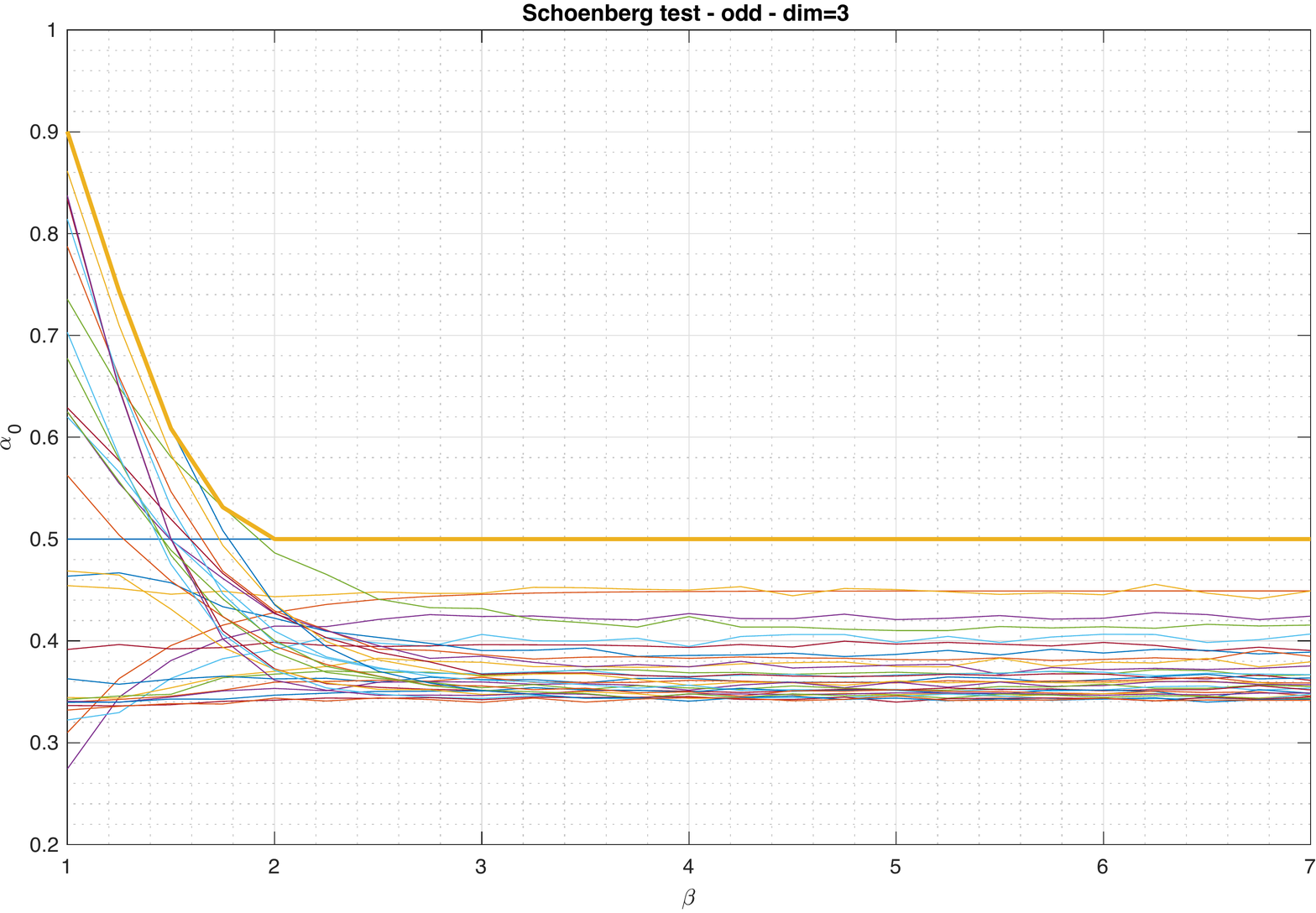}
\caption{{\smallish Fig.2 reproduced by Schoenberg's theorem. }}
\label{fig:alpha}
\end{figure}

\end{document}